\documentclass{article}[14pt]
\usepackage{amsmath}
\usepackage{amssymb}

\newtheorem{prob}{Problem}
\newtheorem{lem}{Lemma}
\newtheorem{thm}{Theorem}
\newtheorem{cor}[lem]{Corollary}

\renewcommand{\le}{\leqslant}
\renewcommand{\ge}{\geqslant}
\renewcommand{\th}{\theta}
\newcommand{\om}{\omega}
\renewcommand{\r}{\rho}

\renewcommand{\k}{\varkappa}
\renewcommand{\l}{\lambda}
\newcommand{\Om}{\Omega}
\newcommand{\se}{\subseteq}
\newcommand{\sd}{\leftthreetimes}
\newcommand{\la}{\langle}
\newcommand{\ra}{\rangle}
\newcommand{\ov}{\overline}
\renewcommand{\a}{\alpha}
\newcommand{\ot}{\otimes}
\newcommand{\be}{\begin{equation}}
\newcommand{\ee}{\end{equation}}
\newcommand{\ve}{\varepsilon}
\newcommand{\ZZ}{\mathbb{Z}}
\newcommand{\NN}{\mathbb{N}}
\newcommand{\FF}{\mathbb{F}}
\newcommand{\ld}{\ldots}
\newcommand\SL{{\rm SL}}
\newcommand\GL{{\rm GL}}
\newcommand\SU{{\rm SU}}
\newcommand\PSL{{\rm PSL}}
\newcommand\PSU{{\rm PSU}}
\newcommand{\ba}[1] {\begin{array}{#1}}
\newcommand{\ea} {\end{array}}

\begin{document}

\large

\begin{center}
{\LARGE \bf Properties of element orders in covers for \\
${\rm L}_n(q)$ and ${\rm U}_n(q)$}

\vspace{0.5cm}

{\sc Andrei\,V.\,Zavarnitsine\footnote{supported by FAPESP, Brazil,
Proc. 06/60766-3; RFBR, Russia, grant 05-01-00797; SB RAS, grant
No.29 for young scientists and Integration Project 2006.1.2}}

\vspace{0.5cm}

{\small  Sobolev Institute of Mathematics, \\
pr. Koptyuga 4, Novosibirsk,  630090, Russia }

\end{center}

\hfill UDC 512.54

\begin{abstract}
We show that if a finite simple group $G$ isomorphic to $\PSL_n(q)$ or $\PSU_n(q)$, where
either $n\ne 4$, or $q$ is prime or even, acts on a vector space over a field of the defining
characteristic of $G$, then the corresponding semidirect product contains an element whose
order is distinct from every element order of $G$. As a consequence, we prove that the group
$\PSL_n(q)$, $n\ne 4$ or $q$ prime or even, is recognizable by spectrum from its covers thus
giving a partial positive answer to Problem 14.60 from the Kourovka notebook.
\end{abstract}

\section {Introduction}

If a group $H$ is a homomorphic image of a finite group  $G$ then
we say that $G$ is a {\it cover} for $H$, or that $G$ {\it covers}
$H$. This paper is devoted to the following problem included in
the Kourovka notebook \cite[Problem 14.60]{k}:

\begin{prob}\label{prob} Suppose that $G$ is a proper
cover for the finite simple group $L={\rm L}_n(q)$, $n\ge 3$. Is it true that $G$ contains an
element whose order is distinct from the order of every element of $L$?
\end{prob}

This problem is related to the recognition of finite groups by
spectrum. Recall that the {\it spectrum} $\om(H)$ of a finite
group $H$ is the set of its elements orders. We call $H$ {\it
recognizable (by spectrum) from its covers} if, for every finite
group $G$ covering $H$, the equality of the spectra
$\om(G)=\om(H)$ implies the isomorphism $G\cong H$. Thus, Problem
\ref{prob} asks if every simple group ${\rm L}_n(q)$, $n\ge 3$, is
recognizable from its covers.

Some special cases of this problem have already been treated elsewhere, e. g. \cite{z,vg,mzz}.
Moreover, the simple groups ${\rm L}_2(q)$ are recognizable from their covers due to
\cite{bs,mz}.

It can be shown (see Lemma \ref{lact}) that the consideration of Problem \ref{prob} may be
reduced to the case where the cover $G$ is  the natural semidirect product $W\sd L$, where $W$
is an elementary abelian $p$-group, $p$ being the defining characteristic for $L={\rm
L}_n(q)$, and the action of $L$ on $W$ is faithful and absolutely irreducible. We prove that
such a $G$ usually contains an element of new order. More precisely, if we denote ${\rm
L}^+_n(q)=\PSL_n(q)$ and ${\rm L}_n^-(q)=\PSU_n(q)$, then the following holds:

\begin{thm}\label{ilu} Let $\ve\in\{+,-\}$
and let $L={\rm L}^{\ve}_n(q)$, $q=p^m$ be a simple group. Suppose that either $n\ge 5$, or
$n=4$ and $q$ is prime, or $n=4$ and $q$ is even. If $L$ acts on a vector space $W$ over a
field of characteristic $p$ then $\om(W\sd L)\ne\om(L)$.
\end{thm}
As follows from the proof, in the case where either $n\ge 5$ or $q$ is an odd prime
we can assert even more: the group $L$ contains a {\it semisimple}
element $g$ of  $p$-maximal order (i.e. such that $p\,|g|\not \in
\omega(L)$) which centralizes a nontrivial vector in $W$.
Moreover, if $n\ge 5$ and $q>3$, such an element $g$ may be chosen
independent of the module $W$. The proof uses the properties of
weights of the irreducible modules for the algebraic group of type
$A_l$.

As a consequence of this result, we have the following (partial)
affirmative solution to Problem~\ref{prob}:

\begin{cor} \label{icor} Let $L={\rm L}_n(q)$
be a simple linear group.
If either $n\ne 4$, or $q$ is prime, or $q$ is even then $L$ is
recognizable by spectrum from its covers.
\end{cor}

Therefore, the only remaining unresolved case for Problem~\ref{prob} is where $L={\rm L}_4(q)$
with $q$ odd and nonprime. We observe that the action of ${\rm L}_4^\ve(q)$ in the defining
characteristic turned out to be a more subtle issue. The above methods will not always work as
there are examples of semidirect products $W\sd L$ which do not contain elements of order $pt$
for $p$-maximal order $t$ coprime with $p$. This means that the action of unipotent elements
of ${\rm L}_4^\ve(q)$ should also be taken account of. For example, in characteristic $p=2$,
let $W$ be the natural module for $L=\SU_4(2)$. Then $\om(W\sd L)\setminus
\om(L)=\{8\}$. There are also more complicated examples of this
kind in odd characteristic.

\section {Preliminaries}

In what follows, we denote by $\FF_q$ a finite field of $q=p^m$
elements. The center of a group $G$ is $Z(G)$.

Let $t>1$ and $n$ be natural numbers and let $\ve\in\{+,-\}$. If
there exists a prime that divides $t^n-(\ve1)^n$ and does not
divide $t^i-(\ve1)^i$ for $1\leqslant i<n$, then we denote this
prime by $t_{[\ve n]}$ and call it a {\it primitive divisor} of
$t^n-(\ve1)^n$. In general, a primitive divisor need neither exist nor
be unique. The following lemma is a generalization of the
well-known Zsigmondy's theorem:

\begin{lem} \label{zhig} Let $t,n>1$ be natural numbers and $\ve\in\{+,-\}$. There exists a primitive divisor $t_{[\ve
n]}$ of $t^n-(\ve1)^n$, except in the following cases:

(i)\, $\ve=+$, $n=6$, $t=2$;

(ii)\, $\ve=+$, $n=2$, and $t=2^l-1$ for some $l\ge 2$;

(iii)\, $\ve=-$, $n=3$, $t=2$;

(iv)\, $\ve=-$, $n=2$, and $t=2^l+1$ for some $l\ge 0$.
\end{lem}
{\it Proof.\/} See \cite[Lemma 5]{mzz}. ~$\blacktriangle$

Let $q$ be a power of a prime and let $\ve\in\{+,-\}$. For
$n\in\NN$, we define the {\it generalized primitive divisor}
$$
q_{[\ve n]}^*=\left\{
\begin{array}{ll}
q_{[\ve n]},& \text{if}\ \ q_{[\ve n]}\  \text{exists}, \\
9, & \text{if}\ \ (\ve,n,q)=(+,6,2), \\
2^l, & \text{if}\ \ (\ve,n,q)=(+,2,2^l-1)\ \text{for}\ l\ge 2, \\
2^l, & \text{if}\ \ (\ve,n,q)=(-,2,2^l+1)\ \text{for}\ l\ge 2. \\
\end{array}
\right.
$$
Observe that $q_{[\ve n]}^*$ is not defined if and only if
\be\label{qnd}
(\ve,n,q)\in\{(-,2,2),(-,2,3),(-,3,2)\}.
\ee
The following assertion follows directly from the above definition:

\begin{lem}\label{div} Suppose that $r=q_{[\ve n]}^*$ is defined. Then
\begin{enumerate}
\item[$(i)$]  $r \ \big|\  \big(q^s-(\ve 1)^s\big )$ \ if
and only if \ \ $n\mid s$;
\item[$(ii)$] for $n>1$,\ we have\ \ $\gcd(r,q-\ve1)=1$, unless $(\ve,n,q)=(+,2,2^l-1)$ or
$(-,2,2^l+1)$;
\item[$(iii)$] if $s\mid n$ and $s>1$ then
\be\label{ex}
r \ \bigm| \frac{q^n-(\ve 1)^n}{q^{n/s}-(\ve 1)^{n/s}},
\ee
unless $(\ve,n,s,q)=(+,6,3,2)$;
\item[$(iv)$] for $n>1$, the group $\SL_n^\ve(q)$ contains an irreducible
element of order $r$.
\end{enumerate}
\end{lem}

In the following lemmas, a quotient of the finite group $\SL_n^\ve(q)$ by a central subgroup
is said to be a group of type $A^\ve_{n-1}(q)$.

\begin{lem}\label{pin} Let $q$ be a power of a prime $p$.
A group of type $A^\ve_{n-1}(q)$
contains an element of order $p^{t+1}$, $t\ge 0$, if and only if
$n\ge p^t+1$.
\end{lem}
{\it Proof.} See \cite[Corollary 0.5]{test} ~$\blacktriangle$

\begin{lem}\label{nom} $(i)$ Let $L={\rm L}_n^\ve(q)$ be a simple group.
Then the numbers
\be\label{nrs}\frac{q^n-(\ve1)^n}{d(q-\ve 1)},\qquad
\frac{q^{n-1}-(\ve1)^{n-1}}{d} \ee are coprime and
and maximal by divisibility elements of $\omega(L)$.

$(ii)$ Let $n\in \NN$, let $q$ be a power of a prime $p$ and let $\ve\in\{+,-\}$. Suppose that
$n=s+b_1+\ld+b_k$, where $s=0$ or $s=p^t$ with $t\ge 0$, $k\ge 0$, and all $b_i$'s are
pairwise coprime and greater than $1$. If $(\ve,q)=(-,2)$ then we assume additionally that
$b_i\ne 2,3$ for all $i$. If $(\ve,q)=(-,3)$ then we assume that $b_i\ne 2$ for all $i$. Put
$r_i=q^*_{[\ve b_i]}$ (which exists due to the restrictions on $b_i$). Then $p^{t+1}r_1\cdot
\ld\cdot r_k\not\in\om(\SL^\ve_n(q))$, where it is assumed that $t=0$ if $s=0$.
\end{lem}
{\it Proof.} We prove item $(ii)$ first. The particular case $s=p^t$ was proved in \cite[Lemma
9]{mzz}; however, we do not exclude it from the consideration as we give here a different
proof. Suppose to the contrary that there is $a\in \SL^\ve_n(q)$ of order $p^{t+1}r_1\cdot
\ld\cdot r_k$. Then $a$ lies in the centralizer $C$ in $\SL^\ve_n(q)$ of the semisimple
element $u=a^{p^{t+1}}$. By \cite[Propositions 7, 8]{cc}, $C$ is a central product of groups
$M_{i,j}$, $i\ge 2$, of type $A^{\ve}_{i-1}(q^{{\mu^{(i)}}_j})$ extended by an abelian group
$T$ of order $\prod_{i,j}(q^{{\mu^{(i)}}_{\!\!j}}-(\ve1)^{{\mu^{(i)}}_{\!\!j}})/(q-\ve 1)$,
where $\mu^{(i)}$ is a partition of $n_i$ and the numbers $n_i$ satisfy $\sum i\, n_i=n$. In
particular, \be\label{smu}
\sum_{i,j}i\,{\mu^{(i)}}_{\!\!j}=n. \ee Observe that $u\in Z(C)$
and $Z(C)$ is an abelian group of order dividing
 $\prod_{i,j}(q^{{\mu^{(i)}}_{\!\!j}}-(\ve1)^{{\mu^{(i)}}_{\!\!j}})$.
Since $|u|=r_1\cdot \ld\cdot r_k$, by Lemma \ref{div}$(i)$ it follows that, for each
$f=1,\ld,k$, there are $i,j$ such that $b_f$ divides ${{\mu^{(i)}}_{\!\!j}}$
\big(if $r_f$ is not prime then we should also recall that $Z(C)$ is in fact a subgroup of the
direct product of cyclic groups of the orders
$q^{{\mu^{(i)}}_{\!\!j}}-(\ve1)^{{\mu^{(i)}}_{\!\!j}}$ for all $i,j$\,\big).  By hypothesis,
all the numbers $b_f$ are coprime and greater than $1$; hence, the sum of those
${{\mu^{(i)}}_{\!\!j}}$ greater than $1$ is at least $b_1+\ld+b_k$.

Now, because $C$ contains an element of order $p^{t+1}$, one of
the above components $M_{i',j'}$ is nontrivial and such that
$i'\ge p^t+1$ by Lemma \ref{pin}. If ${{\mu^{(i')}}_{\!\!j'}}=1$,
we have \be\label{cont} \sum_{i,j}i\,{{\mu^{(i)}}_{\!\!j}}\ge
b_1+\ld+b_k+p^t+1>n, \ee contrary to (\ref{smu}). If
${{\mu^{(i')}}_{\!\!j'}}>1$, then $i'\,{{\mu^{(i')}}_{\!\!j'}}\ge
p^t+1+{{\mu^{(i')}}_{\!\!j'}}$ and we still have (\ref{cont}).
This contradiction completes the proof.

We now prove $(i)$. First, following the above argument we show that the numbers in
(\ref{nrs}) are in $\mu(M)$. Let $k$ denote either of these numbers. It is well-known that $L$
contains an element of order $k$. We show the maximality of $k$ in $\om(L)$ (with respect to
divisibility). We may assume that $(\ve,n,q)\ne(-,3,3),\allowbreak(-,4,2)$, because, for the
groups ${\rm U}_3(3)$ and ${\rm U}_4(2)$, the assertion is readily verified. Suppose that
there is an element $\ov{a}\in L$ of order a multiple of $k$. Then the preimage $a\in
S=\SL_n^\ve(q)$ of $\ov{a}$ lies in the centralizer $C$ in $S$ of a semisimple element $u$ of
order $k$. Observe that the generalized primitive divisor $r=q^*_{[\ve n]}$ (respectively,
$r=q^*_{[\ve (n-1)]}$) exists, since $L\ne {\rm U}_3(3),{\rm U}_4(2)$. Then $u\in Z(C)$ and we
conclude as above that $n$ (respectively, $n-1$) divides some ${\mu^{(i)}}_{\!\!j}$. But then
(\ref{smu}) implies that the decomposition for $n$ is $n=n_1={\mu^{(1)}}_{\!1}$ (respectively,
$n=n_1={\mu^{(1)}}_{\!1}+{\mu^{(1)}}_{\!2}$, where ${\mu^{(1)}}_{\!1}=n-1$ and
${\mu^{(1)}}_{\!2}=1$). In particular, $C$ coincides with its toral part $T$ isomorphic to the
cyclic group of order $(q^n-(\ve1)^n)/(q-\ve 1)$ (respectively, $q^{n-1}-(\ve1)^{n-1}$). Due
to the conjugacy of the maximal tori, $T$ contains the center $Z(S)$ of order $d$ and hence
the order of $\ov{a}\in T/Z(S)$ does not exceed $k$.

We finally show that the numbers in (\ref{nrs}) are coprime.
Denote
$$
x=\frac{(\ve q)^n-1}{\ve q-1},\qquad y=\frac{\ve q-1}{d},\qquad
z=\frac{(\ve q)^{n-1}-1}{\ve q-1}.
$$
Then up to sign the numbers in (\ref{nrs}) coincide with $x/d$ and
$yz$, respectively. Note that $\gcd(x,z)=1$, since
$$\gcd\big((\ve q)^n-1,(\ve q)^{n-1}-1\big)=(\ve q)^{\gcd(n,n-1)}-1=\ve q-1.$$
Also, setting $f(t)=t^{n-1}+\ld+t+1\in \ZZ[t]$ we can find
$g(t)\in \ZZ[t]$ such that $f(t)=(t-1)g(t)+n$. Then the
substitution $t=\ve q$ gives $x=f(\ve q)\equiv n\pmod{\ve q-1}$.
Thus, $\gcd(x,\ve q-1)=\gcd(n,\ve q-1)=d$ and so $\gcd(x/d,y)=1$.
The claim follows from these remarks. ~$\blacktriangle$

\begin{lem}\label{prim}$(i)$ For every
real $x\ge 22$, the interval
$(x/3,x/2]$ contains at least one prime.

$(ii)$ For every real $x\ge 57$, the interval $(2x/3,x-16]$
contains at least one prime.

$(iii)$ For every real $x>45$, the interval $(3x/4,x-8)$ contains
at least one prime.

$(iv)$ For every real $x>27$, the interval $(x/2,x-8)$ contains at
least two primes.

\end{lem}
{\it Proof.} $(i)$ For $x<72$ the assertion is readily verified.
Suppose that $x\ge 72$. There exists $\a\in [0,3)$ such that
$x/3+\a=3a$ for an integer $a>1$. By \cite{ha}, the interval
$(3a,4a)$ contains a prime. It remains to show that $4a\le x/2$.
We have $4a=4(x/3+\a)/3<4x/9+4=x/2-(x/18-4)\le x/2$, since
$x/18-4\ge 0$. Hence, $(i)$ holds.

Items $(i)$--$(iv)$ can be proved in a similar manner, except that
in  $(ii)$ we should use a stronger result that, for every natural
$n\ge 119$, the interval $[n, 1.073n]$ contains at least one
prime, see \cite{rw}. ~$\blacktriangle$

\begin{lem}\label{decex}  $(i)$ For every natural $n\ge 5$
there exists a decomposition
$n=n_1+\ld+n_k$, where $n_1,\ld,n_k$ are pairwise coprime natural
numbers at most one of which is equal to $1$, such that the
following property holds: for every $1\le j\le n$, there is a
decomposition $j=j_1+\ld+j_{k'}$ with $k'\le k$ and an injection
$\eta:\{j_1,\ld,j_{k'}\}\to\{n_1,\ld,n_k\}$ satisfying, for all
$i=1,\ld,k'$, the conditions
\begin{enumerate}
\item[$(a)$] $j_i\le \eta(j_i)$;
\item[$(b)$] if  $\eta(j_i)>1$ then $\gcd(j_i,\eta(j_i))>1$.
\end{enumerate}

$(ii)$ For every natural $n\ge 5$ and $1\le j\le n$ such that $(n,j)\not
\in\{(6,3),(8,3),(8,5)\}$, there exist a decomposition $n=n_1+\ld+n_k$, where $n_1,\ld,n_k$ are
pairwise coprime natural numbers distinct from $2,3$ at most one of which is equal to $1$, and
a decomposition $j=j_1+\ld+j_{k'}$ having the same properties as in $(i)$ and satisfying the
additional requirement that $\gcd(j_i,\eta(j_i))\ne 3$ whenever $\eta(j_i)=6$, where
$i=1,\ld,k'$.
\end{lem}
{\it Proof.} $(i)$ For a natural $m>1$, we denote by $\k(m)$ the largest prime divisor of $m$.

We will show that a stronger fact holds: there exists a required decomposition $n=n_1+\ld+n_k$
with the additional property that \be\label{kap} \k(n_1\ld n_k)\le (n+1)/2.
\ee
We proceed by induction on $n$. Suppose that $n\le 20$.

If $n=5$ then  $n=1+4$ is the required decomposition. Indeed, for $j=1,2,4$, we may set $k'=1$
and $j_1=j$, while, for $j=3$ or $5$, we take $k'=2$ and $j=1+2$ or $j=1+4$, respectively. If
$n=6$ then we decompose $n=1+2+3$. For all $j$, the corresponding decomposition
$j=j_1+\ld+j_{k'}$ is obvious. If $n=7$, we set $n=1+6$. For $j=1,2,3,6$, we set $k'=1$ and,
for $j=4,5,7$, we set $k'=2$ and $j=1+3$, $1+4$, and $1+6$, respectively. If $n=9$, we set
$n=1+8$. For $j=1$ or $j$ even, we set $k'=1$ and, for $j>1$ odd, we set $k'=2$ and $j_1=1$,
$j_2=j-1$.

In all the above cases $n=5,6,7,9$, the injection $\eta$ is
straightforward and the property (\ref{kap}) holds.

Now, for $n=8,10,11,\ld,20$ we define recursively $n=[n-r]+r$ for the respective values
$r=3,5,5,5,7,7,\ld,7,5,3$, where $[n-r]$ denotes the decomposition for $n-r$ already defined.
It is directly verified that all $n_i$ are pairwise coprime, at most one of them is $1$, and
that (\ref{kap}) holds. If $j\le n-r$ then the decomposition $j=j_1+\ld+j_{k'}$ is defined  as
for $n-r$ with the same embedding $\eta$, while if $j>n-r$ then we set $j=[j-r]+r$ and extend
$\eta$ by setting $\eta(r)=r$. (Note that, for $n=13$ and $j=7$, the decomposition $[j-r]$ is
considered empty.)

Suppose that $n\ge 21$. By Lemma \ref{prim}$(i)$, there exists a
prime $r$ such that $(n+1)/3<r\le (n+1)/2$. Since $n-r\ge
(n-1)/2\ge 5$, by induction there exists a decomposition
$n-r=n_1+\ld+n_{k_0}$ satisfying the hypothesis and, additionally,
such that $\k(n_1\ld n_{k_0})\le (n-r+1)/2$. We show that
$n=n_1+\ld+n_{k_0}+r$ is the required decomposition. Since
$r>(n-r+1)/2$, it follows that each prime divisor of each $n_i$ is
less than $r$; in particular, the numbers $n_1,\ld,n_{k_0},r$ are
pairwise coprime, at most one of them is $1$, and $\k(n_1\ld
n_{k_0}r)=r\le (n+1)/2$.

Now let $1\le j \le n$. As above, if $j\le n-r$ then the decomposition $j=j_1+\ld+j_{k'}$ is defined by
induction with the same embedding $\eta$. Suppose that $j\ge n-r+1$. We have $n-r+1\ge r$. If $j=r$ then
we take this equality for a (trivial) decomposition of $j$ and set $\eta(r)=r$; otherwise, $j>r$, i.e.
$1\le j-r
\le n-r$, and we set $j=[j-r]+r$, where the decomposition $[j-r]$ and the embedding $\eta$ on the components
of $[j-r]$ are defined by induction, and set $\eta(r)=r$. From the construction it is clear that all the
requirements on $\eta$ are satisfied, which completes the proof of $(i)$.

We now prove item $(ii)$. In this case, we will find the needed decompositions for $n$ and $j$
such that $\eta(j_i)=n_i$ for all $1=1,\ld,k'$ (we may fix the order of
the summands $n_i$, because $j$ is now fixed.)

$(a)$\ First, we assume that either $j=2$, $3$, $n-2$, $n-3$, or
$(n,j)\in\{(8,4),(16,6),(16,10)\}$. Then the decompositions for $n$ and $j$ are shown in Table
\ref{decpart} in the columns labeled $[n]$ and $[j]$, respectively. Note that, in each case,
due to the restrictions on $n$ and $j$ the components $n_i$ in $n=n_1+\ld+n_k$ are distinct
from $2,3$, and $(n_i,j_i)\ne 3$ whenever $n_i=6$. All the remaining requirements are readily
verified.

\begin{table}[htb]
\caption{A decomposition for  $j=2,3,n-2,n-3$ or
$(n,j)\in\{(8,4),(16,6),(16,10)\}$ \label{decpart}}
{\footnotesize
\begin{center}
\begin{tabular}{l|c|c|c|c|c}
\hline
$(n,j)$ & restrictions on $n$ & $k$ & $[n]$& $k'$& $[j]$  \\
\hline
\hline
$(n,2)$ &
\begin{tabular}{c}
$n\equiv 1\pmod 2$\\
$n\equiv 0\pmod 2$
\end{tabular}
&
\begin{tabular}{c}
$2$\\
$1$
\end{tabular}
&
\begin{tabular}{c}
$(n-1)+1$\\
$n$
\end{tabular}
&
\begin{tabular}{c}
$1$\\
$1$
\end{tabular}
&
\begin{tabular}{c}
$2$\\
$2$
\end{tabular}
\\
\hline
$(n,n-2)$ &
\begin{tabular}{c}
$n\equiv 1\pmod 2$\\$n\equiv 0\pmod 2$
\end{tabular}
&
\begin{tabular}{c}
$2$\\
$1$
\end{tabular}
&
\begin{tabular}{c}
$1+(n-1)$\\$n$
\end{tabular}
&
\begin{tabular}{c}
$2$\\
$1$
\end{tabular}
&
\begin{tabular}{c}
$1+(n-3)$\\$n-2$
\end{tabular}
\\
\hline
$(n,3)$ &
\begin{tabular}{c}
$n\equiv 1\pmod 2$\\$n\equiv 0\pmod 6$\\$n\equiv 2\pmod 6$\\$n\equiv 4\pmod 6$
\end{tabular}
&
\begin{tabular}{c}
$2$\\
$1$\\
$3$\\
$2$
\end{tabular}
&
\begin{tabular}{c}
$1+(n-1)$\\$n$\\$1+4+(n-5)$\\$(n-1)+1$
\end{tabular}
&
\begin{tabular}{c}
$2$\\
$1$\\
$2$\\
$1$
\end{tabular}
&
\begin{tabular}{c}
$1+2$\\$3$\\$1+2$\\$3$
\end{tabular}
\\
\hline
$(n,n-3)$ &
\begin{tabular}{c}
$n\equiv 1\pmod 2$\\$n\equiv 0\pmod 6$\\$n\equiv 2\pmod 6$\\$n\equiv 4\pmod 6$
\end{tabular}
&
\begin{tabular}{c}
$2$\\
$1$\\
$3$\\
$2$
\end{tabular}
&
\begin{tabular}{c}
$(n-1)+1$\\$n$\\$4+(n-5)+1$\\$1+(n-1)$
\end{tabular}
&
\begin{tabular}{c}
$1$\\
$1$\\
$2$\\
$2$
\end{tabular}
&
\begin{tabular}{c}
$n-3$\\$n-3$\\$2+(n-5)$\\$1+(n-4)$
\end{tabular}
\\ \hline
$(8,4)$, $(16,6)$, $(16,10)$ & --- & $1$ & $n$ & $1$ & $j$
\\ \hline
\end{tabular}
\end{center}
}
\end{table}

$(b)$\ We may now assume that $j\ne 2,3,n-2,n-3$ and
$(n,j)\not\in\{(8,4),(16,6),(16,10)\}$. We show by induction on
$n$ that in this case a stronger fact holds: there are needed
decompositions $n=n_1+\ld+n_k$ and $j=j_1+\ld+j_{k'}$ with the
additional requirement that $j_i=n_i$ for all $i=1,\ld, k'$.

Observe that we may also assume that $j\le n/2$. (Otherwise, we
take $n-j$ instead of $j$ and decompose $n=n_1+\ld+n_k$ and
$n-j=n_1+\ld+n_{k'}$. Then $j=n_{k'+1}+\ld+n_{k}$ is the required
decomposition for $j$, while $n$ retains the same decompositions
with the appropriate permutation of the summands.) We consider
three subcases.

$(b.1)$\ Induction base. If $n\le 112$ then it can be directly
verified that, for each admissible pair $(n,j)$, the required
decompositions for $n$ and $j$ have one of the forms
$$
\ba{l}
1)\ \ n=(j)+(n-j), \quad j=(j);\\
2)\ \ n=(1)+(j-1)+(n-j), \quad j=(1)+(j-1);\\
3)\ \ n=(j)+(1)+(n-j-1), \quad j=(j);\\
4)\ \ n=[j]+[n-j], \quad j=[j];\\
\ea
$$
where in 1) --- 3) a summand in parentheses denotes a single
component of the decomposition, while in 4) a summand in brackets
is decomposed as shown in Table \ref{baddec}. For example, if
$(n,j)=(21,7)$ then $\gcd(j,n-j-1)=\gcd(7,13)=1$ and so the
decompositions in 3) satisfy the requirements.

\begin{table}[htb]
\caption{An exceptional decomposition for admissible pairs
$(n,j)$, with $2j\le n\le 112$\label{baddec}} {\footnotesize
\begin{center}
\begin{tabular}{l|l|l||l|l|l||l|l|l}
\hline
$(n,j)$ & $[n-j]$ & $[j]$&  $(n,j)$ & $[n-j]$ & $[j]$&$(n,j)$ & $[n-j]$ & $[j]$       \\
\hline
$(21,6)$& $4+11$&$1+5$    &$(76,36)$& $40$&$7+29$    &$(99,22)$& $77$&$5+17$         \\
$(25,10)$& $4+11$&$1+9$   &$(78,22)$& $56$&$5+17$    &$(99,36)$& $63$&$5+31$          \\
$(34,12)$& $22$&$5+7$     &$(81,6)$& $4+71$&$1+5$    &$(100,12)$&$88$&$5+7$         \\
$(36,15)$& $21$&$4+11$    &$(81,15)$& $7+59$&$15$    &$(100,22)$&$78$&$5+17$        \\
$(45,12)$& $33$&$5+7$     &$(81,36)$& $45$&$7+29$    &$(100,45)$&$55$&$4+41$        \\
$(46,6)$& $11+29$&$6$     &$(85,10)$& $4+71$&$1+9$   &$(105,14)$&$91$&$5+9$           \\
$(46,10)$& $7+29$&$10$    &$(85,15)$& $11+59$& $15$  &$(105,39)$&$5+61$&$39$          \\
$(49,21)$& $5+23$&$21$    &$(85,34)$& $51$&$5+29$    &$(105,40)$&$65$&$7+33$          \\
$(51,6)$& $4+41$&$1+5$    &$(85,35)$& $9+41$&$35$    &$(106,6)$& $11+89$& $6$        \\
$(51,15)$& $7+29$&$15$    &$(85,40)$& $45$&$11+29$   &$(106,10)$&$7+89$&$10$          \\
$(52,18)$& $34$&$5+13$    &$(88,30)$& $58$&$7+23$    &$(106,28)$&$78$&$5+23$        \\
$(55,10)$& $4+41$&$1+9$   &$(91,21)$& $11+59$& $21$  &$(106,36)$&$70$&$13+23$         \\
$(55,15)$& $11+29$& $15$  &$(91,26)$& $65$&$7+19$    &$(106,40)$&$66$&$17+23$           \\
$(55,22)$& $33$&$5+17$    &$(91,28)$& $63$&$5+23$    &$(106,50)$&$56$&$9+41$        \\
$(57,21)$& $5+31$&$21$    &$(91,35)$& $9+47$&$35$    &$(111,6)$& $4+101$& $1+5$       \\
$(64,28)$& $36$&$5+23$    &$(91,39)$& $5+47$&$39$    &$(111,12)$&$99$&$5+7$         \\
$(66,26)$& $40$&$7+19$    &$(92,14)$& $5+73$&$14$    &$(111,15)$&$7+89$&$15$        \\
$(69,18)$& $51$&$5+13$    &$(93,24)$& $69$&$5+19$    &$(111,33)$&$5+73$&$33$        \\
$(70,24)$& $46$&$5+19$    &$(96,20)$& $76$&$7+13$    &$(111,36)$&$75$&$7+29$        \\
$(76,6)$& $11+59$&$6$     &$(99,21)$& $5+73$&$21$    &$(111,45)$&$7+59$&$45$          \\
\cline{4-9}
$(76,10)$& $7+59$&$10$\\
\cline{1-3}
\end{tabular}
\end{center}
}
\end{table}

$(b.2)$\ Suppose that $n\ge 113$ and $j\le 2n/5$. Then $n-j\ge 3n/5>57$ and by
Lemma \ref{prim}$(ii)$ there is a prime $r$ such that
\be\label{rcnd}
2(n-j)/3<r\le n-j-16.
\ee
Clearly, $r\ne 2,3$. Moreover, we have
\be\label{jcnd}
j\le 2(n-j)/3<r,
\ee
and hence
\be\label{rg}
(n-j)/2< 2(n-j)/3<r.
\ee
Now the pair $(n-r,j)$ satisfies the induction hypothesis. Indeed, by (\ref{rcnd}) we have
$j\ne 2,3,n-r-2,n-r-3$,  and  $n-r\ge j+16>16$. Hence, $(n-r,j)$ is an admissible
pair and by induction we have $n-r=n_1+\ld+n_k$ and $j=n_1+\ld+n_{k'}$, where $k'\le k$
and $n_i$'s are
pairwise coprime, distinct from $2,3$, and at most one of them being $1$. Since
$r>j$ by (\ref{jcnd}), and $r>n-r-j$ by (\ref{rg}), it follows that $r$ is greater than,
hence coprime with, all $n_i$'s. Thus $n=n_1+\ld+n_k+r$ and $j=n_1+\ld+n_{k'}$
are the required decompositions for $n$ and $j$.

$(b.3)$\  Suppose that $n\ge 113$ and $j> 2n/5$. Then $j>45$ and
by Lemma \ref{prim}$(iii)$ there is a prime $s$ such that
\be\label{scnd} 3j/4<s<j-8. \ee Since $(n-j)/2<3j/4$, we have
\be\label{j8g} (n-j)/2<s. \ee Due to $j\le n/2$ we also have
$(n-j)/2\ge n/4>27$ and so Lemma \ref{prim}$(iv)$ implies that
there is a prime $r$ distinct from $s$ such that \be\label{rg2}
(n-j)/2<r<n-j-8. \ee Consider the pair $(n-s-r,j-s)$. By
(\ref{rg2}) and (\ref{scnd}), we have $j-s>8>2,3$; $n-j-r>8>2,3$;
and $n-s-r=(n-j-r)+(j-s)>8+8=16$. Thus $(n-s-r,j-s)$ is an
admissible pair and by induction we have $n-s-r=n_1+\ld+n_k$  and
$j-s=n_1+\ld+n_{k'}$. By (\ref{j8g}) and (\ref{rg2}), we have
$s,r>(n-j)/2\ge j/2$. Thus $s,r$ are distinct from, hence coprime
with, all $n_i$'s. Therefore, $n=s+n_1+\ld+n_k+r$ and
$j=s+n_1+\ld+n_{k'}$ are the required decompositions. This
completes the proof of the lemma. ~$\blacktriangle$

We emphasize that the difference of case $(ii)$ from case $(i)$ of Lemma \ref{decex} is
not only in the requirement that $n_i\ne 2,3$ but also in that the decomposition
for $n$ depends on the number $j$.

Let $r$ be a prime, $G$ a finite group, and $g\in G$. We say that
$g$ is an element of {\em $r$-maximal order} if $r\,|g|\not
\in\om(G)$. Examples of $r$-maximal orders for the groups $\SL_n^\ve(q)$ are given in
Lemma \ref{nom}$(ii)$

\begin{lem}\label{elem} Let $S=\SL_n^\ve(q)$, with $q=p^m$.

$(i)$ Let $n\ge 5$ and $q>3$. Then  $S$ contains a semisimple element $g$ of $p$-maximal order
such that $\la g \ra\cap Z(S)=1$ and, for every $0\le j\le n$, the product of some $j$
distinct characteristic values of $g$ (in the natural $n$-dimensional representation) equals
$1$.

$(ii)$ Let $n\ge 4$, $0\le j\le n$, and $(\ve,n,q)\ne (-,4,2)$. Then $S$ contains a semisimple element
$g$ of $p$-maximal order such that $\la g \ra\cap Z(S)=1$ and the product of some $j$ distinct
characteristic values of $g$ (in the natural $n$-dimensional representation) equals $1$.

\end{lem}
{\it Proof.} $(i)$ Let $n=n_1+\ld+n_k$ be the decomposition whose existence is stated in Lemma
\ref{decex}$(i)$. Then $S$ includes a naturally embedded subgroup isomorphic to
$\SL_{n_1}^\ve(q)\times\ld\times\SL_{n_k}^\ve(q)$.  By Lemma~\ref{div}$(iv)$ and in view of
the restriction $q>3$, we may choose an element $g_i\in\SL_{n_i}^\ve(q)$ of order
\be\label{ri}
r_i=\left\{ \ba{rl}
1,& n_i=1; \\
q_{[\ve n_i]}^*,& n_i>1. \ea\right.
\ee
Set $g=g_1\ld g_k\in S$ (so that $g$ is the direct sum of diagonal blocks $g_i$).
By the coprimality of $n_1,\ld,n_k$, we have
$|g|=r_1\ld r_k$, and hence $|g|$ is $p$-maximal
by Lemma \ref{nom}. Observe also that by Lemma \ref{div}$(ii)$
either $|g|$ is coprime with $q-\ve1$, or $q=2^l\pm 1$ and there is
$1\le i_0\le k$ such that $n_{i_0}=2$. However, in the latter case, we must
have $k\ge 2$, since $n\ge 5$. These remarks imply that $\la g \ra\cap Z(S)=1$
by the construction of $g$.

Clearly, the set of characteristic
values for $g$ is the union of those for $g_i$ which have the form
\be\label{chval} \{\th_i^{\strut },\th_i^{\ve q},\th_i^{(\ve
q)^2},\ld,\th_i^{(\ve q)^{n_i-1}}\}, \ee for some $\th_i\in F^\times$ of order $r_i$,
$i=1,\ld,k$, where $F$ is the algebraic closure of $\FF_p$.

If $j=0$, we set $g$ to be any semisimple element  of $p$-maximal order such that $\la g \ra\cap
Z(S)=1$. Let $1\le j \le n$ and $j=j_1+\ld+j_{k'}$ as stated in Lemma
\ref{decex}$(i)$. Without loss of generality (renumbering, if
necessary, the summands in $n=n_1+\ld+n_k$) we may assume that $\eta(j_i)=n_i$, $i=1,\ld,k'$,
where $\eta$ is defined in Lemma \ref{decex}$(i)$. It is sufficient to show that the product
of some distinct $j_i$ values in (\ref{chval}) equals $1$. We may assume that $n_i>1$
(otherwise, $\th_i=1$ and the claim holds). Then $d_i=\gcd(j_i,n_i)>1$ by the property $(b)$
in Lemma \ref{decex}$(i)$. Observe that, by Lemma~\ref{div}$(iii)$, we have
\be
\label{rid}
r_i \ \bigm|
\frac{(\ve q)^{n_i}-1}{(\ve q)^{n_i/d_i}-1}=1+x+x^2+ \ld +x^{d_i-1},
\quad x=(\ve q)^{n_i/d_i}.
\ee
In particular, the set (\ref{chval}) is the union of $f=n_i/d_i$
mutually disjoint subsets
$$
\{\th_i^{\vphantom{x}}\,,\th_i^x,\ld,\th_i^{x^{d_i-1}}\},\ \{\th_i^{\ve
q},\th_i^{x(\ve q)},\ld,\th_i^{x^{d_i-1}(\ve q)}\},\
\ld\ ,\{\th_i^{(\ve q)^{f-1}},\th_i^{x(\ve
q)^{f-1}},\ld,\th_i^{x^{d_i-1}(\ve q)^{f-1}}\},
$$
in each of which the product of all elements  equals $1$ due to
(\ref{rid}). Since $j_i/d_i\le f$ by the property $(a)$ of Lemma
\ref{decex}$(i)$, it follows that the union of arbitrary $j_i/d_i$
of the above subsets gives the required set of $j_i$ distinct characteristic values whose
product equals $1$.

$(ii)$ As above, we may assume that $j>0$. First, suppose that $n\ge 5$ and
$(n,j)\not\in\{(6,3),(8,3),\allowbreak (8,5)\}$. Then we decompose $n=n_1+\ld+n_k$ and
$j=j_1+\ld+j_{k'}$ as stated in Lemma \ref{decex}$(ii)$. As above, there exists an element
$g=g_1\ld g_k\in S$ of $p$-maximal order $r_1\ld r_k$, where the $r_i$'s defined by (\ref{ri})
exist due to the restrictions $n_i\ne 2,3$. We now repeat the rest of the argument of $(i)$ to
show that there are $j$ distinct characteristic values of $g$ whose product equals $1$.

If $(n,j)\in\{(8,3),(8,5)\}$ and $(\ve,q)\ne (-,2)$ then because of (\ref{qnd}) we may allow a
summand of $n$ to equal $3$. So we decompose $n=5+3$, $j=j$ (trivial decomposition). If
$(n,j,\ve,q)\ne (6,3,+,2)$ then the divisibility (\ref{rid}) holds by (\ref{ex}) and we may
repeat the above argument.

If $n=4$ and $(\ve,q)\ne (-,2)$ then again $3$ is a possible summand for $n$ and so we
decompose $n=1+3$ if $j=1$ or $3$, and decompose trivially $n=4$ if $j=2$ or $4$ and proceed
as above.

Suppose that $S=\SL_6(2)$ and $j=3$. Then $S$ contains an element $g$ of $2$-maximal order
$21$ whose characteristic values are $\nu_i=\th^{2^{i-1}}$, $i=1,\ld,6$, where $\th\in
\ov{\FF}_2^\times$ is of order $21$. Observe that $Z(S)=1$ and
the product of $3$ characteristic values $\nu_1,\nu_3,\nu_5$ of $g$ equals
$\th\th^4\th^{16}=1$, as required.

Finally, let $S=\SU_8(2)$ and $j=3$ or $5$. Then $S$ contains an element $g$ of order $45$
(which is $2$-maximal) and whose block-diagonal form $g=g_1g_2g_3$ has three blocks of sizes
$4$, $3$, $1$ and the characteristic values $\nu_1,\ldots,\nu_8$ of $g$ are \be\label{gch}
\underbrace{\strut\th^3,\ (\th^3)^{-2},\ (\th^3)^4,\ (\th^3)^{-8}}_{g_1},\
\underbrace{\strut\th^{-5},\ (\th^{-5})^{-2},\ (\th^{-5})^4}_{g_2},\
\underbrace{\strut\th^{30}}_{g_3}\ ; \ee
where $\th\in \ov{\FF}_2^\times$ is of order $45$. Then
$\nu_1\nu_3\nu_8=\nu_2\nu_4\nu_5\nu_6\nu_7=1$ and so there are $j$ characteristic values of
$g$ whose product equals $1$. Since $Z(S)=1$, the claim follows. ~$\blacktriangle$

\begin{lem} \label{frob_action} If a Frobenius group
$KC$ with kernel $K$ and cyclic complement $C=\langle c \rangle$
of order $n$ acts faithfully on a vector space  $V$ over a field
of nonzero characteristic $p$ coprime with the order of $K$, then
the minimal polynomial of  $c$ on $V$ is $x^n-1$. In particular,
the semidirect product $V\sd C$ contains an element of order
$p\cdot n$ and $\mbox{dim} C_V(c)>0$.
\end{lem}

{\it Proof.} See  \cite[Lemma 1]{m94}. ~$\blacktriangle$

\begin{lem}\label{absi} A group $H$
 is recognizable by spectrum from its covers
if and only if $\om(H)=\om(G)$ for every split extension $G=N\sd
H$, where $N$ is an elementary abelian $r$-group for some $r$ and
$H$ acts on $N$ absolutely irreducibly.
\end{lem}
{\it Proof.} Let $G$ be a proper cover for $H$ of minimal order
such that $\om(H)=\om(G)$. By \cite[Lemma 12]{zl}, we may assume
that $G=N\sd H$, where $H$ acts on the elementary abelian
$r$-group $N$ irreducibly. Suppose that this action is not
absolutely irreducible. Let $F$ be a finite splitting field for
$H$ of characteristic $r$ and consider a proper submodule $N_0$ of
the reducible $FH$-module $N\ot _{\FF_p}F$. It is sufficient to
show that $\om(N_0\sd H)=\om(H)$. Suppose to the contrary that
$n_0h\in N_0\sd H$ is an element of order not belonging to
$\om(H)$. Then the element $1+h+h^2+\ld+h^{|h|-1}$ considered as a
linear transformation of $N_0$ is nonzero. But then it is also
nonzero as a linear transformation of $N$ and hence $G$ contains
an element $nh$ of order $|n_0h|$, a contradiction.
~$\blacktriangle$

\begin{lem} \label{non-f} Let $r$ be a prime
and let  $L={\rm L}_n^\ve(q)$ be a simple group, where $q=p^m$ and $\ve\in\{+,-\}$. Then
$\om(\ZZ_r\times L)\not\subseteq\om(L)$.
\end{lem}
{\it Proof.} Denote by $a_1$ and $a_2$ the numbers in (\ref{nrs}).
By Lemma \ref{nom}$(i)$, there is $i=1,2$ such that $r\nmid a_i$
and $ra_i\not\in\om(L)$. Since $ra_i\in \om(\ZZ_r\times L)$, the
claim follows. ~$\blacktriangle$

\begin{lem}\label{lact}  Let $L={\rm L}_n(q)$ be a simple linear group, where $q=p^m$.
Then $L$ is recognizable by spectrum from its covers if and only if $\om(L)=\om(G)$ for every
split extension $G=N\sd L$, where $N$ is an elementary abelian $p$-group and $L$ acts on $N$
faithfully and absolutely irreducibly.
\end{lem}
{\it Proof.} By Lemma \ref{absi}, we may assume that $G=N\sd L$ where $N$ is an elementary
abelian $r$-group for some $r$ and $L$ acts on $N$ absolutely irreducibly. By Lemma
\ref{non-f}, $L$ acts faithfully. The image in $L$ of the parabolic subgroup of $\SL_n(q)$ of
the shape $q^{n-1}\!:\!\GL_{n-1}(q)$ contains by \cite[Lemma 5]{vg} a Frobenius subgroup $KC$
with elementary abelian kernel $K$ of order $q^{n-1}$ and cyclic complement $C$ of order
$a=(q^{n-1}-1)/d$, where $d=\gcd(n,q-1)$. If $r\ne p$ then, by Lemma \ref{frob_action}, we
have $ra\in
\om(G)$. However, $ra\not\in\om(L)$ by Lemma \ref{nom}$(i)$.
Hence, $r=p$. ~$\blacktriangle$

\section{Weights of irreducible $\SL_n(F)$-modules}

In this section, we recall some facts from the representation theory of algebraic groups. For
details, see e. g. \cite{jan}.

Let $G=\SL_n(F)$, where $F$ is an algebraically closed field of characteristic $p$. Then $G$
is a simply connected simple algebraic group of type $A_l$, where $l=n-1$. Denote by $\om_0$
the zero weight and by $\om_1,\ld,\om_l$ the fundamental weights of $G$ (with respect to a
fixed maximal torus of $G$ and a system of positive roots). Let
$\Om=\ZZ\om_1\oplus\ld\oplus\ZZ\om_l$ be the weight lattice and $\Delta$ the root system of
$G$ with the set $\a_1,\ldots,\a_l$ of simple roots. Also, let
$\Om_0=\ZZ\a_1\oplus\ld\oplus\ZZ\a_l$ be the set of {\it radical weights} and
$\Om^+=\{a_1\om_1+\ld+a_l\om_l\in \Om\ |\ a_1\ge 0,\ld,a_l\ge 0\}$ be the set of {\it dominant
weights}. The weights in the set $\Om^+_k=\{a_1\om_1+\ld+a_l\om_l\in \Om\ |\ 0\le
a_1<k,\ld,0\le a_l< k\}$ are called {\it $k$-restricted}, where $k$ is usually a power of $p$.

For an irreducible (rational, finite dimensional) $G$-module $L$,
denote by $\Om(L)$ the set of weights of $L$ and by $\l(L)$ the
highest weight of $L$. It is known that $\l(L)\in\Om^+$ and each
dominant weight is the highest weight of some irreducible module
$L$. The irreducible $G$-module of highest weight $\l$ is
customarily denoted by $L(\l)$. Obviously,
$\Om(L(\om_0))=\{\om_0\}$. The module $L$ is called {\it
$p$-restricted} if $\l(L)\in \Om_p^+$. The modules $L(\om_i)$,
$i=1,\ld,l$ are called the {\it microweight modules}. The
structure of the microweight modules is well known and
described in the following lemma (see, e.g. \cite[II.2.15]{jan}):

\begin{lem}\label{mwm} Let $G=\SL_n(F)$ and let $V=F^n$
be the natural $G$-module with the canonical basis $e_1,\ld,e_n$. Choose the diagonal subgroup
$H$ for a fixed maximal torus of $G$. Then $e_i$ is an eigenvector for $H$ with the
corresponding weight $\ve_i$. Choose a system of positive roots $\{\ve_i-\ve_j\ |\ 1\le i<j\le
n \}$. Then, for $1\le k< n$, we have $\om_k=\ve_1+\ld+\ve_k$ and the microweight module
$L(\om_k)$ is isomorphic to the $k$-th exterior power $\wedge^k V$ and has the set of weights
$$\Om(L(\om_k))=\{\ve_{i_1}+\ve_{i_2}+\ld+\ve_{i_k}\ |\ 1\le i_1<i_2<\ld<i_k\le n \}.$$
\end{lem}

The following assertion is a refinement of \cite[Proposition 2.3]{zs} for groups of type
$A_l$.

\begin{lem}\label{zal} Let $G=\SL_n(F)$ and let $L$ be an irreducible
$p$-restricted $G$-module. Write
$\l(L)=a_1\om_1+\ld+a_l\om_l$. Suppose that $i\in\{0,1,\ldots,l\}$
is the uniquely defined integer such that \be\label{iai}
a_1+2a_2+\ld+la_l\equiv i \pmod{l+1}. \ee Then
$\Om(L(\om_i))\se\Om(L)$.
\end{lem}
{\it Proof.}  By \cite[Proposition II.2.4]{jan}, $\Om(L)$ lies in a single coset of
$\Om:\Om_0$ and, by \cite[Proposition 2.3]{zs}, $\Om(L)$ contains $\om_0$ if $\l(L)$ is
radical, and includes $\Om(L(\om_i))$ for some $i=1,\ld,l$, otherwise. In the latter case, the
index $i$ is uniquely defined, since the weights $\om_0,\om_1,\ld,\om_l$  form a transversal
of $\Om:\Om_0$ by
\cite[VIII, \S 7.3, Proposition 8]{bour}.

Therefore, it remains to observe that an arbitrary weight
$\l=a_1\om_1+\ld+a_l\om_l \in \Om$ lies in the coset $\om_i+\Om_0$
if and only if (\ref{iai}) holds. Indeed, if $\l=\a_i$ is a simple
root then $(a_1,\ld,a_l)$ is the $i$th row of the Cartan matrix of
type $A_l$ and (\ref{iai}) is directly verified. By the above,
every weight $\l$ is in $\om_i+\Om_0$ for some $i$ and adding or
subtracting positive roots from $\l$ preserves the relation
(\ref{iai}) and the coset $\om_i+\Om_0$. The claim follows from
these remarks. ~$\blacktriangle$

\section{Proofs of the main results}

We first give a proof of Theorem \ref{ilu}.

{\it Proof.} By Lemma \ref{absi}, we may assume that $W$ is
absolutely irreducible as a module for $L$. Moreover, by \cite[Theorem 43]{s}, we may
assume that $W$  is a restriction to $S=\SL_n^\ve(q)$ of an
irreducible module $L(\l)$ for $G=\SL_n(F)$, where $\l\in\Om_q^+$,
and $F=\ov{\FF}_p$. Set $l=n-1$.

$(a)$\ We first consider the case where $n\ge 5$ and $q>3$. It is sufficient to show that
there is a semisimple element $g\in S$ of $p$-maximal order
such that $\la g \ra \cap Z(S)=1$
which centralizes a nonzero vector $w\in W$. Indeed, if this is the case then
the element $w\ov{g}\in W\sd L$ has order $p\,|g|\not\in\om(L)$,
due to $\om(L)\se\om(S)$, where $\ov{g}$ is the image of $g$ in $L$.

We choose for a fixed maximal torus of $G$ the diagonal subgroup $H$. Write
$$\l=\l_0+p\l_1+\ld+p^{m-1}\l_{m-1},\qquad \l_i\in\Om^+_p.$$
By the Steinberg Tensor Product theorem \cite[Theorem 41]{s}, we have
\be\label{ste} L(\l)\cong
L(\l_0)\ot L(\l_1)^\r\ot\ld\ot L(\l_{m-1})^{\r^{m-1}}, \ee where $\r$ denotes the twisting by
the Frobenius map corresponding to the automorphism $x\mapsto x^p$ of $F$. By Lemma \ref{zal},
for each $i=0,\ld,m-1$, there exists $k_i\in \{0,\ld,l\}$ such that
$\Om(L(\om_{k_i}))\se\Om(L(\l_i))$. In particular, the set $\Om(L(\l))$ contains all possible
weights of the form
\be\label{pw} \mu_0+p\mu_1+\ld
+p^{m-1}\mu_{m-1}, \qquad \mu_i\in \Om(L(\om_{k_i})). \ee By Lemma
\ref{mwm}, $\mu_i$ can be an arbitrary sum of $k_i$ distinct
weights in $\{\ve_1,\ld,\ve_n\}$.

Let $g\in S$ be the semisimple element whose existence is stated in Lemma \ref{elem}$(i)$. Then there is
$a\in G$ such that $h={}^ag\in H$. By Lemma \ref{elem}$(i)$, there are $k_i$ distinct characteristic
values of $g$ whose product equals $1$ and we define $\mu_i$ to be the sum of the corresponding $k_i$
weights $\ve_j$ so that $\mu_i(h)=1$ for all $i=0,\ld,m-1$. Denote by $\mu$ the sum  (\ref{pw}) with
$\mu_i$ just defined. Then we have $\mu\in \Om(L(\l))$ and
$$
\mu(h)=\mu_0(h)\mu_1(h)^p\ld\mu_{m-1}(h)^{p^{m-1}}=1.
$$
Let $w_0\in W$ be a weight vector for $G$ of weight $\mu$ so that
$w_0h=\mu(h)w_0=w_0$. Set $w=w_0a$. Then
$$
wg=w_0ag=w_0ha=w_0a=w.
$$
Thus $g$ is the required semisimple element of $S$.

$(b)$\ We now suppose that $n\ge 4$, $q$ is prime, and $(\ve,n,q)\ne (-,4,2)$. In this case, $\l\in
\Om_q^+=\Om_p^+$. By Lemma~\ref{zal}, there exists $j\in \{0,\ld,l\}$ such that
$\Om(L(\om_j))\se\Om(L(\l))$. Hence, by Lemma~\ref{mwm}, $\Om(L(\l))$ contains the sum of
arbitrary $j$ distinct weights in $\{\ve_1,\ld,\ve_n\}$. Now,  we choose by Lemma
\ref{elem}$(ii)$ a semisimple element $g\in S$ of $p$-maximal order such that $\la g \ra\cap
Z(S)=1$ and the product of some $j$ distinct characteristic values of $g$ equals $1$. There is
$a\in G$ such that $h={}^ag\in H$ and so the product of some $j$ characteristic values of $h$
equals $1$ as well. We set $\mu$ equal to the sum of the corresponding $j$ weights
$\{\ve_1,\ld,\ve_n\}$ so that $\mu(h)=1$. (Then $\mu
\in\Om(L(\l))$ by the above.) If now $w_0\in W$ is a weight vector
for $G$ then, as in case $(a)$, $wg=w$, where $w=w_0a$ and so $g$
is as required.

We emphasize that, in this cases, the principal difference from
case $(a)$ is that the module $W$ is $p$-restricted and that the
choice of $g$ depends on $W$.

$(c)$\ Finally, let $n=4$ and let $q$ be even. By Lemma 6 in \cite{mzz}, $L$
contains a Frobenius subgroup $KC$ with kernel $K$ of order
$q_{[\ve 4]}$ and cyclic complement $C$ of order $4$. By Lemma
\ref{non-f}, $KC$ acts faithfully on $W$ and hence Lemma
\ref{frob_action} implies that  $2|C|\in \om(W\sd L)$. However,
$2|C|=8\not\in\om(L)$ by Lemma \ref{pin}. This completes the proof of the theorem.
~$\blacktriangle$

Corollary \ref{icor} is now a direct consequence of Lemma \ref{lact}
and Theorem \ref{ilu}.

{\em Acknowledgement.\/} The author is thankful to D. O. Revin for reading the manuscript of
the paper and making a number of valuable remarks.


\begin{thebibliography}{100}
\bibitem{k} Unsolved problems in group
theory, {\it The Kourovka notebook}, 14th ed., Sobolev Inst. Mat.
(Novosibirsk), 1999.

\bibitem{z} A. V. Zavarnitsine, Weights of the irreducible
$\SL_3(q)$-modules in defining characteristic. (Russian) {\em
Sibirsk. Mat. Zh.}, {\bf 45}, N\,2 (2004), 319--328. English
translation in {\em Sib. Math. J.}, {\bf 45}, N\,2 (2004),
261--268

\bibitem{vg} A. V. Vasil'ev and M. A. Grechkoseeva,
On Recognition by Spectrum of Finite Simple Linear Groups over
Fields of Characteristic 2. (Russian) {\em Sibirsk. Mat. Zh.},
{\bf 46}, N\,4 (2005), 749--758. English translation in {\em Sib.
Math. J.}, {\bf 46}, N\,4 (2005), 593--600.

\bibitem{mzz}  V. D. Mazurov, A. V. Zavarnitsine, On element orders
in coverings of the simple groups ${\rm L}_n(q)$ and ${\rm U}_n(q)$, {\em Proceedings of the
Steklov Institute of Mathematics}, Suppl. 1 (2007), 145--154.


\bibitem{bs} R. Brandl, W. Shi,
The characterization of $\PSL_2(q)$ by its element orders, {\em J.
Algebra}, {\bf 163} (1994), N\,1, 109--114.

\bibitem{mz} V. D. Mazurov, A. V. Zavarnitsine,
Element orders in coverings of symmetric and alternating groups.
(Russian) {\em Algebra i Logika}, {\bf 38}, N\,3 (1999), 296--315.
Translation in {\em Algebra and Logic}, {\bf 38}, N\,3 (1999),
159--170.

\bibitem{cc} R. W. Carter, Centralizers of semisimple elements in
the finite classical group, {\em Proc. London Math. Soc.} (3),
{\bf 42}, N\,1 (1981), 1--41

\bibitem{test} D. M.  Testerman, $A_1$-type overgroups of elements of order $p$ in
semisimple algebraic groups and the associated finite groups, {\em
J. Algebra}, {\bf 177}, N\,1,  (1995), 34--76.

\bibitem{zs} I. D. Suprunenko, A. E. Zalesskii, Fixed vectors for
elements in modules for algebraic groups, {\em Internat. J. Algebra Comput.}, {\bf 17},
N\,5--6 (2007), 1249--1261.

\bibitem {zl} A. V. Zavarnitsine,  Recognition of the simple groups
${\rm L}_3(q)$ by element orders. {\em J. Group Theory}, {\bf 7}, N\,1, (2004), 81--97.

\bibitem{bour} N. Bourbaki, \'El\'ements de math\'ematique.
Fasc. XXXVIII: Groupes et alg\`ebres de Lie. Chapitre VII:
Sous-alg\`ebres de Cartan, \'el\'ements r\'eguliers. Chapitre
VIII: Alg\`ebres de Lie semi-simples d\'eploy\'ees. Actualit\'es
Scientifiques et Industrielles, No. 1364. Hermann, Paris, (1975).

\bibitem{jan} J. C. Jantzen, Representations of algebraic groups.
Second edition. Mathematical Surveys and Monographs, 107. American
Mathematical Society, Providence, RI, 2003.

\bibitem{ha} D. Hanson,
On a theorem of Sylvester and Schur, {\em Canad. Math. Bull.} {\bf
16} (1973), 195--199.

\bibitem {s} R. Steinberg,
Lectures on Chevalley groups. (Russian) Biblioteka sbornika
"matematika", Moscow: Verlag "Mir", (1975).

\bibitem{rw} H.Rohrbach, J.Weis, Zum finiten Fall des Bertrandschen
Postulates. J. Reine Angew. Math., 214/215 (1964), 432--440.

\bibitem{m94}  
V. D. Mazurov, On the set of orders of elements of a finite group. (Russian)
{\em Algebra i Logika}, {\bf 33}, N\,1 (1994), 81--89. Translation in {\em
Algebra and Logic} {\bf 33}, N\,1 (1994), 49--55.

\end{thebibliography}
\end{document}